\documentclass[a4paper,abstracton]{scrartcl}

\usepackage{authblk}

\usepackage{mathtools}
\usepackage{amssymb}
\usepackage{bm}

\usepackage{algpseudocode}
\usepackage{algorithm}
\algnewcommand{\Initialise}[1]{%
	\State \textbf{Initialise:}\space\parbox[t]{.8\linewidth}{\raggedright #1}
}

\usepackage{natbib}
\usepackage{hyperref}
\usepackage{xcolor}
\definecolor{darkgreen}{RGB}{0,127,173}
\definecolor{darkblue}{RGB}{0,0,200}
\hypersetup{colorlinks, linkcolor={darkblue}, citecolor={darkgreen}, urlcolor={darkblue}}

\usepackage{graphicx}
\graphicspath{{figures/}}  % path to directory containing all figures
\usepackage{caption}  % allows subfigures
\usepackage{subcaption}  % allows subfigures

\allowdisplaybreaks  % allows page breaks within equations
\usepackage{multirow}  % table formatting
\usepackage{verbatim}  % allows blocks to be commented out

\newcommand{\M}{{\mathcal{M}}}
\newcommand{\U}{{\mathcal{U}}}
\newcommand{\X}{{\mathcal{X}}}
\begin{document}

\title{A faster exact method for solving the robust multi-mode resource-constrained project scheduling problem}

\author[1]{Matthew Bold\footnote{Corresponding author, email: m.bold1@lancaster.ac.uk}}
\author[2]{Marc Goerigk}
\affil[1]{STOR-i Centre for Doctoral Training, Lancaster University, Lancaster, UK}	
\affil[2]{Network and Data Science Management, University of Siegen, Siegen, Germany}	

\date{}

\maketitle

\begin{abstract}
	This paper presents a mixed-integer linear programming formulation for the multi-mode resource-constrained project scheduling problem with uncertain activity durations. We consider a two-stage robust optimisation approach and find solutions that minimise the worst-case project makespan, whilst assuming that activity durations lie in a budgeted uncertainty set. Computational experiments show that this easy-to-implement formulation is many times faster than the current state-of-the-art solution approach for this problem, whilst solving over 40\% more instances to optimality over the same benchmarking set.
\end{abstract}

\noindent\textbf{Keywords:} project scheduling; optimisation under uncertainty; robust optimisation; budgeted uncertainty

%%%%%%%%%%%%%%%%%%%%%%%%%%%%%%%%%%%%%%%%%%%%%%%%%%%%%%%%%%%%%%%%%%%%%%%%%

\section{Introduction}
The multi-mode resource-constrained project scheduling problem (MRCPSP) is a generalisation of the widely-studied resource-constrained project scheduling problem (RCPSP) to include multiple processing modes for each activity. The inclusion of these modes allows the modelling of situations in which there is more than one way of executing project activities, with each option having its own duration and resource requirements. The MRCPSP consists of selecting the processing modes and start times for a given set of activities, subject to a set of precedence constraints and limited resource availability, with the objective of minimising the overall project duration, known as the makespan.

In this paper we consider the MRCPSP under uncertain activity durations and model it using a two-stage adjustable robust optimisation framework. In this setting, a first-stage problem is solved to determine activity mode selections and make activity sequencing decisions to resolve resource conflicts. Following this, the actual activity durations are realised and a complete schedule is computed. The aim of the two-stage adjustable robust MRCPSP is to find a feasible first-stage solution (i.e. mode selection and sequencing decisions) in order to minimise the realised worse-case makespan, as computed in the second-stage. We refer to this problem as the robust MRCPSP.

A number of papers in recent years have applied this two-stage robust optimisation approach to the RCPSP. First to use this approach were \cite{artigues2013robust}, who presented an iterative scenario-relaxation algorithm for this problem with the objective of minimising the worst-case absolute regret. More recently, \cite{bruni2017adjustable} introduced a Benders'-style decomposition approach for solving the robust RCPSP with the objective of minimising the worst-case project makespan. This work was extended in \cite{bruni2018computational}, which presented a computational study comparing an additional Benders' decomposition approach against a primal decomposition algorithm. Most recently, \cite{bold2021compact} introduced a compact reformulation of the robust RCPSP and presented results which showed the superior computational performance of that formulation over the iterative decomposition-based methods developed in the two preceding papers. 

The application of this two-stage robust optimisation approach for the MRCPSP is a very recent development. To the best of our knowledge, the only existing paper to consider this problem is \cite{balouka2021robust}, in which the Benders' decomposition approach introduced by \cite{bruni2017adjustable} for the robust RCPSP has been extended for application to the MRCPSP. Mirroring that extension, this paper adapts the compact formulation developed by \cite{bold2021compact} to the MRCPSP, with the aim of achieving similarly superior computational performance over the Benders' solution approach. 

Following a formal description of the two-stage robust MRCPSP in Section \ref{section:problem_description}, we outline the proposed compact formulation for this problem in Section \ref{section:compact_formulation} and present a strengthened version of the Benders' decomposition solution approach from \cite{balouka2021robust} in Section \ref{section:benders}. Results from a computational comparison of these two approaches are detailed in Section \ref{section:experiments}.

\section{Problem description}\label{section:problem_description}

\begin{figure}[h]
	\centering
	\includegraphics[scale=0.9]{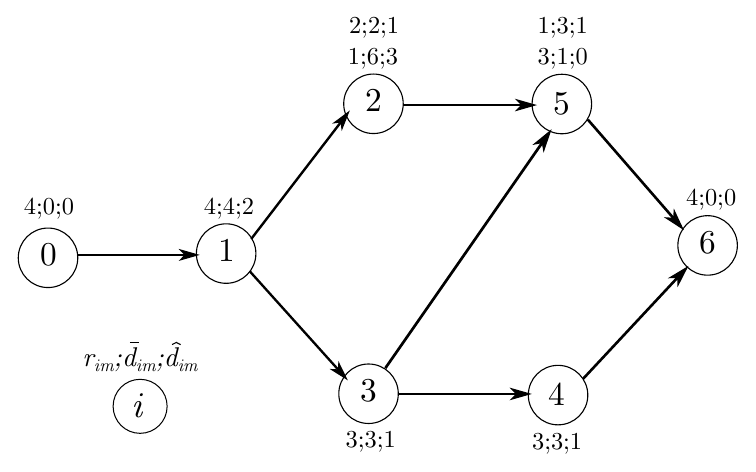}
	\caption{An example instance involving five non-dummy activities and a single renewable resource with availability $R_1=4$. Activities 2 and 5 each have two available processing modes, whilst the other activities have only a single mode.}
	\label{fig:mrcpsp_network}
\end{figure}

A project consists of a set of non-preemptive activities $V=\{0,1,\dots,n,n+1\}$, where 0 and $n+1$ denote the dummy source and sink activities respectively. Each activity $i\in V$ has a set of available processing modes given by $M_i=\{1,\dots,|M_i|\}$. The nominal duration of activity $i$ when executed in mode $m\in M_i$ is given by $\bar{d}_{im}$, whilst its worst-case duration is given by $\bar{d}_{im}+\hat{d}_{im}$, where $\hat{d}_{im}$ is its maximum durational deviation. Each mode for activity $i$, $m\in M_i$, has an associated renewable resource requirement of $r_{imk}$ for each $k\in K$, where $K$ is the set of renewable resource types involved in the project. Each renewable resource $k\in K$ has an availability of $R_{k}$ at each time period in the project horizon. As well as renewable resource requirements, mode $m$ of activity $i$ also has a non-renewable resource requirement of $r^{'}_{imk}$ for each non-renewable resource type $k\in K^{'}$, with each non-renewable resource having an overall availability of $R^{'}_k$ for the entire project horizon. Additionally, the project is subject to a set of strict finish-to-start precedence constraints given by $E$, where $(i,j)\in E$ enforces that activity $i$ must finish before activity $j$ can begin. These form a project network that can be represented using a directed graph $G(V,E)$. Figure \ref{fig:mrcpsp_network} shows an example instance involving five non-dummy activities and a single renewable resource.

For a given choice of processing modes for each activity $\bm{m}=(m_1,\dots,m_n)$, we assume that the activity durations lie somewhere in a budgeted uncertainty set of the form
$$\U_{\bm{m}}(\Gamma)=\Bigg\{\bm{d}_{\bm{m}}\in\mathbb{R}_+^{|V|}:d_{im_i}=\bar{d}_{im_i}+\xi_i\hat{d}_{im_i},\,0\leq \xi_i \leq 1,\,\forall i\in V,\,\sum_{i\in V}\xi_i \leq \Gamma\Bigg\}.$$
Introduced by \cite{bertsimas2004price}, the motivation of the budgeted uncertainty set is to control the pessimism of the solution by introducing a robustness parameter $\Gamma$ to limit the number of jobs that can simultaneously achieve their worst-case durations. Observe that when $\Gamma=0$, each activity takes its nominal duration and the resulting problem is the deterministic MRCPSP. At the other extreme, when $\Gamma=n$, every activity takes its worst-case duration, in which case the uncertainty set becomes equivalent to an interval uncertainty set. In this case the problem can again be solved as a deterministic MRCPSP instance considering only worst-case durations.

For a given set of activity modes $\bm{m}$, a forbidden set is defined to be any subset of activities $F_{\bm{m}}\subseteq V$ that are not precedence-related, such that $\sum_{i\in F_{\bm{m}}}r_{im_ik}>R_k$ for at least one resource $k\in K$. That is, a forbidden set is a collection of activities that cannot be executed in parallel only because of resource limitations. Applying the main representation theorem of \cite{bartusch1988scheduling}, for a particular choice of activity modes $\bm{m}$, a solution to the MRCPSP can be defined by a set of additional precedences $X_{\bm{m}}\subseteq V^2\setminus E$ such that the extended precedence network $G(V,E\cup X_{\bm{m}})$ is acyclic and contains no forbidden sets. Such an extension to the project network is referred to as a \textit{sufficient selection}.

The aim of the two-stage robust MRCPSP is therefore to determine activity modes and a corresponding sufficient selection in order to minimise the worst-case project makespan. For the case where activity durations lie in a budgeted uncertainty set that we consider in this paper, this problem can be written as
\begin{align}
	\min_{\bm{m}\in \M,\,X_{\bm{m}}\in \X_{\bm{m}}} \max_{\bm{d}\in\U_{\bm{m}}(\Gamma)}\ & \min S_{n+1}\label{eqn:robust_mrcpsp_1}\\
	\text{s.t. } & S_0 = 0\\
	& S_j - S_i \geq d_{im_i} \qquad \forall (i,j)\in E\cup X_{\bm{m}}\\
	& S_i \geq 0 \qquad \forall i\in V,\label{eqn:robust_mrcpsp_4}
\end{align}
where $\M\subseteq \mathbb{N}^n$ represents the set of all possible combinations of activity processing mode selections, and $\X_{\bm{m}}$ is the set of all possible sufficient selections for the choice of processing modes given by $\bm{m}$. 

\begin{figure}[htb]
	\makebox[\textwidth][c]{
	\begin{subfigure}[t]{.5\textwidth}
		\centering
		\includegraphics[width=0.95\textwidth]{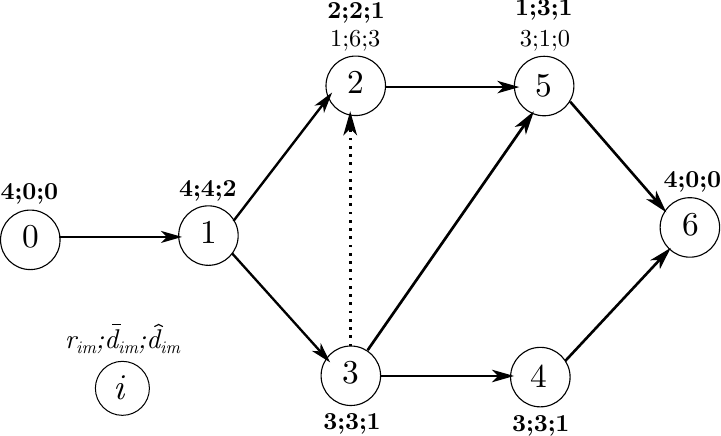}
		\caption{Activity mode choices are highlighted in bold and sufficient selection is given by dashed arc $\{(3,2)\}$.}\label{fig:mrcpsp_network_solution}
	\end{subfigure}%
	\hspace{15mm}
	\begin{subfigure}[t]{.6\textwidth}
		\centering
		\includegraphics[width=1.00\textwidth]{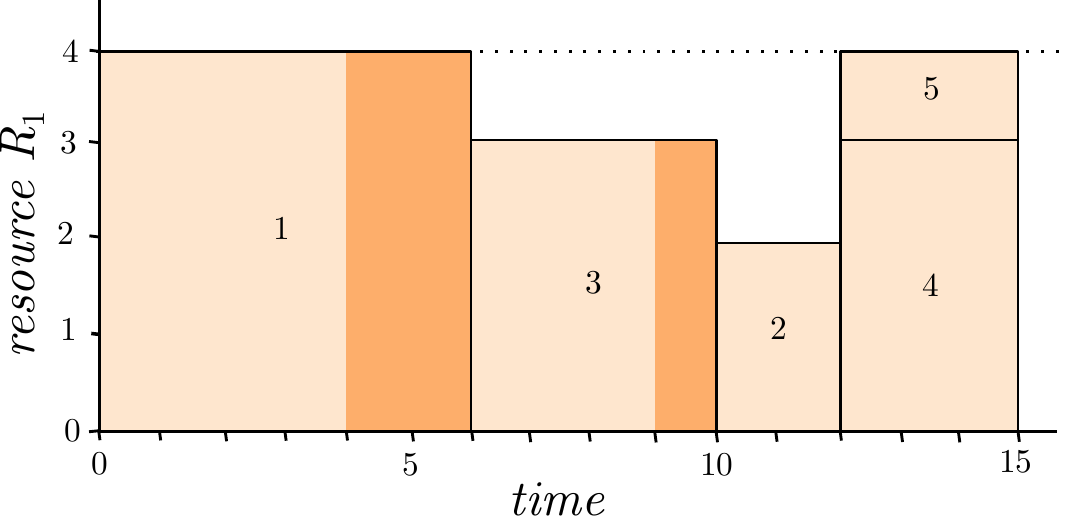}
		\caption{Worst-case schedule corresponding to optimal robust solution, where activities 1 and 3 have been delayed.}\label{fig:mrcpsp_schedule}
	\end{subfigure}
	}
	\caption{Optimal solution to the example instance from Figure \ref{fig:mrcpsp_network}.}
	\label{fig:mrcpsp_solution}
\end{figure}

An optimal solution to the instance shown in Figure \ref{fig:mrcpsp_network} is given by Figure \ref{fig:mrcpsp_network_solution}, where the mode choices are shown in bold, and the sufficient selection is given by $\{(3,2)\}$. The worst-case schedule for this optimal solution is shown in Figure \ref{fig:mrcpsp_schedule}, where activities 1 and 3 have been delayed to achieve a worst-case makespan of 15 (note that delaying a combination of activity 1 and any other activity results in the same worst-case makespan).

\section{A compact formulation}\label{section:compact_formulation}

We propose solving the robust MRCPSP using an extended version of the mixed-integer programming formulation developed by \cite{bold2021compact} for solving the two-stage adjustable robust RCPSP. This formulation combines the first and second-stage problems into a single compact formulation and was shown to significantly outperform the strongest iterative decomposition-based methods for that problem.

This formulation is constructed by recasting the adversarial subproblem of determining the worst-case activity durations for a given set of mode choices and sufficient selection as a longest path problem on an augmented project network. This directed acyclic graph is formed of $\Gamma+1$ connected copies of the original project network, with each level being used to account for a delay to a single activity in the project (a similar construction has also been applied in \cite{bendotti2019anchor}). Having reformulated the adversarial subproblem as a longest-path problem, strong duality can be employed to enable its insertion it into a formulation for the first-stage problem, and ultimately arrive at a complete compact formulation for the full problem. Given that the derivation of this formulation is more or less identical to the derivation of the corresponding formulation for the robust RCPSP, we omit it here and instead refer the reader to \cite{bold2021compact}.

In the resulting compact formulation, $y_{ij}$ variables are used to define a complete set of transitive precedences between the project activities that are formed by the original project precedence constraints and the additional precedences introduced to resolve resource conflicts, i.e. $\{(i,j),\, i,j\in V:y_{ij}=1\}=T(E\cup X_{\bm{m}})$, where $T(\cdot)$ is used to denote the transitive closure of a set. Variables $x_{im}$ determine the processing mode for each activity, whilst resource flow variables $f_{ijk}$ track the amount of renewable resource $k$ that is transferred from activity $i$ to activity $j$ upon its completion. The $S_{i\gamma}$ variables are dual variables associated with the longest-path subproblem, which in this formulation are used to track the start time of activities under the worst-case activity durations. Using these variables, the compact formulation for the robust MRCPSP can be written as   
\begin{align}
\min \ & S_{n+1,\Gamma}\label{eqn:compact1}\\
\text{s.t. } & S_{00}=0\label{eqn:compact2}\\
	     & S_{j\gamma} - S_{i\gamma} \geq \bar{d}_{im}x_{im} - N(1-y_{ij}) \qquad \forall (i,j)\in V^2,\,\forall m\in M_i,\,\gamma=0,\dots,\Gamma\label{eqn:compact_bigM1}\\
	     & S_{j,\gamma+1} - S_{i\gamma} \geq (\bar{d}_{im} + \hat{d}_{im})x_{im} - N(1-y_{ij})\nonumber \\
	     & & \hspace{-8cm} \forall (i,j)\in V^2,\, m\in M_i,\,\gamma = 0,\dots,\Gamma-1\label{eqn:compact_bigM2}\\
	     & y_{ij}=1 \qquad \forall (i,j) \in E\cup\{(n+1,n+1)\}\label{eqn:compact5}\\
	     & y_{ij} + y_{ji} \leq 1 \qquad \forall i,j\in V,\, i<j\label{eqn:compact_trans1}\\
	     & y_{ij} = y_{ip} + y_{p j} - 1 \qquad \forall i,j,p\in V,\,i\neq j\neq p\label{eqn:compact_trans2}\\
	     & f_{ijk}\leq P_{ijk}y_{ij} \qquad \forall (i,j)\in V^2,\,i\neq n+1,\,j\neq0,\, \forall k \in K\label{eqn:compact6}\\
	     & \sum_{i\in V\setminus\{n+1\}}f_{ijk}=\sum_{m\in M_j}r_{jmk}x_{jm} \qquad \forall j \in V\setminus\{0\} ,\, \forall k\in K\label{eqn:compact7}\\
	     & \sum_{j\in V\setminus\{0\}}f_{ijk}=\sum_{m\in M_i}r_{imk}x_{im} \qquad \forall i \in V\setminus\{n+1\},\,\forall k\in K\label{eqn:compact8}\\
	     & \sum_{m\in M_i}x_{im} = 1 \qquad \forall i\in V\label{eqn:compact9}\\
	     & \sum_{m\in M_i}r^{'}_{imk}x_{im} \leq R^{'}_k \qquad \forall i \in V,\,k\in K^{'}\label{eqn:compact_nonrenewable}\\
	     & S_{i\gamma}\geq 0 \qquad \forall i\in V,\,\gamma \in 0,\dots, \Gamma\label{eqn:compact10}\\
	     & y_{ij}\in\{0,1\} \qquad \forall (i,j)\in V^2\label{eqn:compact11}\\
	     & f_{ijk}\geq 0 \qquad \forall (i,j)\in V^2 ,\, \forall k \in K\label{eqn:compact12}\\
	     & x_{im}\in\{0,1\} \qquad \forall i\in V,\,m\in M_i\label{eqn:compact13},
\end{align}
where $N=\sum_{i\in V}\max_{m\in M_i}(\bar{d}_{im}+\hat{d}_{im})$ is an upper bound on the minimum makespan, and $P_{ijk}=\min\{\max_{m\in M_i}r_{imk},\,\max_{m\in M_j} r_{jmk}\}$ is the maximum possible flow of resource $k$ from $i$ into $j$.

Constraints (\ref{eqn:compact2})-(\ref{eqn:compact_bigM2}) are the dual constraints corresponding to the longest-path variables in the adversarial subproblem, and serve to ensure that the activity start times respect the project precedences as well as any delays to activity durations. Constraints (\ref{eqn:compact5}) capture the original project precedence relationships. Although not required for the correctness of the model, constraints (\ref{eqn:compact_trans1}) and (\ref{eqn:compact_trans2}) are additional transitivity constraints that have been included because they were shown in \cite{bold2021compact} to provide significant improvements to the computational performance of the model. Constraints (\ref{eqn:compact6})-(\ref{eqn:compact8}) are resource flow constraints which ensure that the transfer of resources between activities follows precedence constraints and that resources are conserved as they flow into and out of each activity in the network. Constraints (\ref{eqn:compact9}) allows only one processing mode to be selected for each activity. Finally, non-renewable resource constraints are enforced by (\ref{eqn:compact_nonrenewable}).

With polynomially many variables and constraints this formulation can be implemented straightforwardly using standard mathematical optimisation software such as Gurobi or CPLEX.

\section{A Benders' decomposition approach}\label{section:benders}

In this section we outline an alternative approach to solving the robust MRCPSP based on Benders' decomposition. A Benders'-type approach for solving the robust MRCPSP was first presented in \cite{balouka2021robust}, and as mentioned previously, this itself was based on the approach for solving the robust RCPSP used in \cite{bruni2017adjustable}.

The main idea behind this approach is to decompose the full problem into its two stages: 1. a master problem that determines activity processing modes and resolves resource conflicts, and 2. a subproblem that takes the solution from the master problem and evaluates it by finding the worst-case makespan for the resulting network by solving a longest path problem. Since the master problem does not account for all the problem uncertainty at once, its solution forms a lower bound to the optimal objective value of the original problem. Meanwhile, the solution to the subproblem is feasible to the original problem, and therefore provides an upper bound. These two problems are solved iteratively, with optimality cuts being added to the master problem each iteration based on the solution to the subproblem until the lower and upper bounds meet.

After initially implementing the Benders' algorithm as it is presented in \cite{balouka2021robust}, it was discovered that its performance could be strengthened by replacing the master problem formulation used in that implementation with a stripped-down version of the compact formulation (\ref{eqn:compact1})-(\ref{eqn:compact13}) in which the uncertainty is removed. In this section we present our strengthened implementation of the Benders' approach for solving the robust MRCPSP. A comparison of these two Benders' implementations is included in the results in Section \ref{section:experiments}. For the details of the original implementation of the Benders' decomposition algorithm for the robust MRCPSP, see \cite{balouka2021robust}.

\subsection{The master problem}

The role of the master problem is to determine a choice of activity processing modes $\bm{m}$, as well as a sufficient selection $X_{\bm{m}}$ to resolve any resulting resource conflicts. This problem is solved without the direct consideration of the uncertain activity durations, and instead, uncertainty is accounted for in the subproblem and communicated back to the master problem with the use of optimality cuts. Hence, to remove the uncertainty from the model, (\ref{eqn:compact1})-(\ref{eqn:compact13}) is modified by considering only a single level of the augmented project network. Using this formulation, the master problem at iteration $t$ can be written as
\begin{align}
\min \ & \eta\label{eqn:master_first}\\
\text{s.t. } & \eta \geq S_{n+1}\\
	     & S_{0}=0\label{eqn:master2}\\
	     & S_{j} - S_{i} \geq \bar{d}_{im}x_{im} - N(1-y_{ij}) \qquad \forall (i,j)\in V^2,\,\forall m\in M_i\label{eqn:master3}\\
	     & \eta \geq \lambda(\bm{x}, \bm{y}, V^{*\ell}) \qquad \forall \ell=1,\dots,t-1\label{eqn:master_cuts}\\
	     & \textnormal{(\ref{eqn:compact5})-(\ref{eqn:compact_nonrenewable}),\,(\ref{eqn:compact11})-(\ref{eqn:compact13})}\\
	     & S_i\geq 0 \qquad \forall i\in V,\label{eqn:master_last}
\end{align}
where (\ref{eqn:master_cuts}) are the optimality cuts generated from the solutions of the resulting subproblems at each of the previous iterations.

The optimal objective value of the master problem (\ref{eqn:master_first})-(\ref{eqn:master_last}), denoted by $\eta^*$, forms a lower bound to the original robust MRCPSP problem. Note that for the first iteration, the solution to the master problem corresponds to an optimal solution to the MRCPSP with nominal activity durations.

\subsection{The subproblem}
 
The subproblem at iteration $t$ computes the worst-case makespan for the mode choices and sufficient selection from the master problem at iteration $t$, denoted by $\bm{m}^{*t}$ and $X_{\bm{m}^{*t}}$ respectively. This is done by solving a longest-path problem through the extended project network $T(E\cup X_{\bm{m}^{*t}})$, in which the activity durations (i.e. arc lengths) can be chosen from the budgeted uncertainty set $\U_{\bm{m}^{*t}}(\Gamma)$. This problem can be formulated as 
\begin{align}
	V^{*t} = &\max \sum_{(i,j)\in T(E\cup X_{\bm{m}^{*t}})} \bar{d}_{im^{*t}_i}\alpha_{ij} + \hat{d}_{im^{*t}_i}w_{ij}\label{eqn:subproblem_first}\\
	& \sum_{(i,n+1)\in T(E\cup X_{\bm{m}^{*t}})} \alpha_{i,n+1} = 1\\
	& \sum_{(0,i)\in T(E\cup X_{\bm{m}^{*t}})} \alpha_{0,i} = 1\\
	& \sum_{(i,j)\in T(E\cup X_{\bm{m}^{*t}})} \alpha_{ij} - \sum_{(j,i)\in T(E\cup X_{\bm{m}^{*t}})} \alpha_{ji} = 0 \qquad \forall i \in V\setminus\{0,n+1\}\\
	& w_{ij} \leq \xi_i \qquad \forall (i,j)\in T(E\cup X_{\bm{m}^{*t}})\\
	& w_{ij} \leq \alpha_{ij} \qquad \forall (i,j)\in T(E\cup X_{\bm{m}^{*t}})\\
	& \sum_{i\in V}\xi_i \leq \Gamma\\
	& \alpha_{ij}\in \{0,1\} \qquad \forall (i,j)\in T(E\cup X_{\bm{m}^{*t}})\\
	& w_{ij}\geq 0 \qquad \forall (i,j)\in T(E\cup X_{\bm{m}^{*t}})\\
	& 0\leq \xi_i \leq 1 \qquad \forall i\in V. \label{eqn:subproblem_last}
\end{align}

Variables $\alpha_{ij}$ define the longest path through the network, which we denote as $\pi^{*t}=\{(i,j),\,i,j\in V:\alpha_{ij}=1\}$, and has length $V^{*t}$. Variables $\xi_i$ are the activity delay variables, whilst $w_{ij}$ are used to linearise the formulation. Note that this formulation of the subproblem is identical to the one implemented by \cite{balouka2021robust}.

\subsection{Optimality cuts} \label{section:optimality_cuts}

Following the solution to the subproblem, a valid cut is generated and added to the master problem for the next iteration. This cut simply forces an alternative solution in the master problem if it is to achieve a better objective value in the next iteration. 

Given the current best lower bound for the original problem, $LB$, the mode selection from the master problem, $\bm{m}^{*t}$, and the optimal objective value of the subproblem, $V^{*t}$, the cut at iteration $t$ can be written as
\begin{align}
	& \eta \geq (V^{*t} - LB)\cdot\sum_{(i,j)\in \pi^{*t}}\bigg(1/3 (y_{ij} + x_{i,m_i^{*t}} + x_{j,m_j^{*t}}) - (3 - y_{ij} - x_{i,m_i^{*t}} - x_{j,m_j^{*t}})\bigg) & \nonumber \\
	& \hspace{8.1cm} - (V^{*t} - LB)\cdot(|\pi^{*t}|-1) + LB.\label{eqn:benders_cut}
\end{align}

\cite{balouka2021robust} show that this constraint is a valid optimality cut for the problem, and that the number of these cuts that need to be added to the master problem before finding an optimal solution is finite.

An overview of the implementation of the Benders' solution approach outlined in this section is presented in Algorithm \ref{alg:benders}.

\begin{algorithm}[h] 
\caption{Benders' decomposition algorithm.}
\begin{algorithmic}[1]
\Initialise{Set $LB=-\infty$, $UB=+\infty$ and $t=1$.}
\While{$UB > LB$}
	\State \textbf{Solve master problem (\ref{eqn:master_first})-(\ref{eqn:master_last})} 
	\State \hspace{\algorithmicindent} Get objective value $\eta^{*t}$, processing modes $\bm{m}^{*t}$ and precedences $\bm{y}_{\bm{m}^{*t}}$
	\State \hspace{\algorithmicindent} If $\eta^{*t}>LB$, update $LB\leftarrow \eta^{*t}$ 
	\State \textbf{Solve subproblem (\ref{eqn:subproblem_first})-(\ref{eqn:subproblem_last})}
	\State \hspace{\algorithmicindent} Get objective value $V^{*t}$ and longest path $\pi^{*t}$
	\State \hspace{\algorithmicindent} If $V^{*t}<UB$, update $UB\leftarrow V^{*t}$
	\State \textbf{Add cut (\ref{eqn:benders_cut}) to master problem}
	\State Update $t\leftarrow t+1$
\EndWhile
\State \Return $UB$
\end{algorithmic}
\label{alg:benders}
\end{algorithm}

\subsection{Example}

To demonstrate its implementation, we use the Benders' decomposition approach to solve the example shown in Figure \ref{fig:mrcpsp_network} with $R_1=4$ and $\Gamma=2$. The algorithm solves this instance in 6 iterations, and the solution information from the master and subproblem at each of these iterations is shown in Table \ref{table:benders_example}.

The first iteration of the master problem solves the nominal instance, assuming no delays to the activity durations. This solution opts to schedule activities 2 and 3 in parallel, which can be achieved by setting $m_2=2$, and the sufficient selection for is given by $\{(4,5)\}$. This first solution is evaluated in the subproblem to have a worst-case makespan of 16. In the second iteration, having added the first optimality cut, the master problem finds the solution shown in Figure \ref{fig:mrcpsp_network_solution}. By calculating the schedule shown in Figure \ref{fig:mrcpsp_schedule}, the subproblem evaluates the objective value of this solution to be 15. Although this solution is optimal, the algorithm requires a further four iterations to prove that this is the case. Note that there is no need to solve the subproblem in the final iteration since after solving the master problem and updating $UB$, we have that $LB=UB$.
\begin{table}[h]
\centering
{\renewcommand{\arraystretch}{1.2}  % vertical padding
\begin{tabular}{llllllll}
	\hline \hline
	$t$ & $LB$ & $UB$ & $\eta^{*t}$ & $V^{*t}$ & \multicolumn{1}{c}{$\bm{m}^{*t}$}     & \multicolumn{1}{c}{$\pi^{*t}$}  \\
	\hline
1         & 11 & 16 & 11          & 16                                  & 1,2,1,1,2 & 0,1,2,5,6 \\
2         & 12 & 15 & 12          & 15                                  & 1,1,1,1,1 & 0,1,3,2,4,6 \\
3         & 12 & 15 & 12          & 15                                  & 1,1,1,1,1 & 0,1,2,3,4,6 \\
4         & 13 & 15 & 13          & 18                                  & 1,2,1,1,1 & 0,1,2,5,6 \\
5         & 13 & 15 & 13          & 16                                  & 1,1,1,2,1 & 0,1,3,4,2,5,6 \\
6         & 15 & 15 & 15          & -                                   & 1,1,1,1,1 & - \\
	\hline \hline
\end{tabular}
\caption{Solution information at each of the six iterations of the Benders' algorithm required to solve example instance in Figure \ref{fig:mrcpsp_network}.}
\label{table:benders_example}}
\end{table}

\section{Computational experiments and results} \label{section:experiments}
 
In this section we compare results from using the compact formulation and the Benders' decomposition approach to solve uncertain MRCPSP instances. A complete set of the raw results used to generate the tables and plots presented here, in addition to the source code used to implement these experiments, can be found at \url{https://github.com/boldm1/robust-mrcpsp}.

The instances used in this computational study have been created from the deterministic $j10$, and $j20$ MRCPSP instances from the PSPLIB (\cite{kolisch1997psplib}, \url{https://www.om-db.wi.tum.de/psplib/}). Note that these are the same instance sets as used for the experiments in \cite{balouka2021robust}. The $j10$ set contains a total of 536 instances each involving 10 activities, and the $j20$ set contains a total of 554 instances each involving 20 activities. We introduce uncertainty into these deterministic instances by setting the maximum durational deviation for each activity to be $\hat{d}_{im}=\lfloor0.7\times \bar{d}_{im}\rfloor$ for each mode $m\in M_i$. These uncertain instances have then been solved using both methods for a range of robustness levels $\Gamma$. In particular we solve the $j10$ instances for $\Gamma\in\{0,3,5,7\}$ and the $j20$ instances for $\Gamma\in\{0,5,10,15\}$. 

Both the compact reformulation and Benders' decomposition approach have been implemented in Python 3.9.2 and solved using Gurobi 9.0.1 running on a single core of a 2.30 GHz Intel Xeon CPU. A time limit of 2 hours per instance was imposed on each of the solution methods.

We begin by observing the effect of the robustness parameter $\Gamma$ on the average optimal objective values shown in Tables \ref{table:j10_objvals} and \ref{table:j20_objvals}. The values in these tables have been computed only over the instances for which an optimal solution was found for all the values of $\Gamma$ (i.e. all 536 instances in the $j10$ set, and 477 instances of the $j20$ set). As we would expect, the solution cost increases in a concave manner as $\Gamma$ increases. 

\begin{table}[h]
	\centering
	{\renewcommand{\arraystretch}{1.2}  % vertical padding
	\begin{subtable}{.45\linewidth}
	\begin{tabular}{lllll}
		\hline\hline
		$\Gamma$  & 0     & 3     & 5     & 7 \\
		\hline
		obj. & 16.84 & 25.34 & 26.35 & 26.46 \\               
		\hline\hline
	\end{tabular}
	\caption{$j10$}
	\label{table:j10_objvals}
	\end{subtable}%
	\qquad
	\begin{subtable}{.45\linewidth}
	\begin{tabular}{lllll}
		\hline\hline
		$\Gamma$  & 0     & 5     & 10     & 15 \\
		\hline
		obj. & 24.61 & 37.89 & 38.69 & 38.69 \\               
		\hline\hline
	\end{tabular}
	\caption{$j20$}
	\label{table:j20_objvals}
	\end{subtable}
	\caption{Average optimal objective values across instances in the $j10$ and $j20$ instance sets for different values of $\Gamma$.}
	\label{table:objvals}}
\end{table}

Table \ref{table:benders_vs_compact} reports the percentage of instances solved to optimality (\textit{\%sol.}), the average percentage optimality gap over instances for which a feasible solution was found (\textit{gap}), and the average solution time in seconds (\textit{sol. time}) for the different choices of $\Gamma$ across the two instance sets, for both the Benders' decomposition approach and the compact formulation. Note that if no optimal solution was found within the time limit by one of the solution methods, the solution time was recorded as the time limit value of 7200 seconds. Additionally, for the Benders' approach, the average number of (completed) iterations (\textit{it.}) and the average time per (completed) iteration (\textit{it. time}) are also reported. The average of the reported values across the different values for $\Gamma$ are also given for each instance set. 

% Explanation of Benders' getting stuck in first iteration
% Note also that when the Benders' algorithm is unable to find a feasible solution to a particular instance within the given time limit, this means that the algorithm was unable to solve the first iteration of the master problem. Hence this first iteration was not completed and is not counted in the reported \textit{it.} and \text{it. time} values. These instances are however recorded as having a solution time of 7200 seconds which is the reason for the disparity between the average iteration time and the average solution time in the $\Gamma=0$ row of the $j20$ set.

Firstly, the results in Table \ref{table:benders_vs_compact} show that for both methods, the instances tend to increase in difficulty as $\Gamma$ increases. We can also see that for the nominal problems, i.e. when $\Gamma=0$, the Benders' approach has a slight edge on the compact formulation. This is what we would expect given the relationship between the Benders' master problem formulation and the compact formulation. However, for all the non-zero values of $\Gamma$ across both instance sets, the compact formulation performs significantly better than the Benders' approach, solving a greater proportion of instances with dramatically reduced computation times. 

\begin{table}[h]
\centering
\small  % font size
{\renewcommand{\arraystretch}{1.2}  % vertical padding
\begin{tabular}{lrrrrrrrrrrrr}
		     \hline \hline
		     &      & & \multicolumn{5}{c}{Benders'}                                       & & \multicolumn{3}{c}{Compact formulation} \\
		     \cline{4-8} \cline{10-12}
		     & $\Gamma$ & & \textit{\%sol.} & \textit{gap} & \textit{its.} & \textit{it. time} & \textit{sol. time} & & \textit{\%sol.}   & \textit{gap}   & \textit{sol. time}   \\
		     \hline
	\multirow{4}{*}{$j10$} & 0    & & 100.0    & 0.00     & 1            & 0.96          & 1.1	    & & 100.0      & 0.00       & 1.8             \\
			       & 3    & & 81.7     & 3.13     & 67           & 3.76          & 1551.4       & & 100.0      & 0.00       & 5.0             \\
			       & 5    & & 79.7     & 4.20     & 76           & 3.89          & 1688.4       & & 100.0      & 0.00       & 4.6             \\
			       & 7    & & 79.5     & 4.52     & 78           & 3.98          & 1701.1       & & 100.0      & 0.00       & 6.5             \\
		     \cline{4-8} \cline{10-12}
			       &      & & 85.2	   & 2.96     & 56	     & 3.15 	     & 1235.5	    & & 100.0      & 0.00	& 4.47 		  \\
		     \hline
	\multirow{4}{*}{$j20$} & 0    & & 89.9     & 0.00     & 1            & 171.08        & 885.5        & & 89.5       & 1.89       & 925.1           \\
			       & 5    & & 50.9     & 10.21    & 220          & 133.32        & 4169.3       & & 88.6       & 1.87       & 1041.5          \\
			       & 10   & & 50.4     & 10.99    & 225          & 138.49        & 4244.6       & & 87.9       & 2.58       & 1118.5          \\
			       & 15   & & 49.5     & 11.21    & 216          & 146.26        & 4300.3       & & 86.3       & 2.91       & 1260.6          \\
		     \cline{4-8} \cline{10-12}
			       &      & & 60.2	   & 8.10     & 165	     & 147.28	     & 3399.9       & & 88.1       & 2.31       & 1086.4          \\
		     \hline \hline
\end{tabular}
}
\caption{Comparison of the Benders' decomposition approach and the compact reformulation across the $j10$ and $j20$ instance sets and for different values of $\Gamma$.}
\label{table:benders_vs_compact}
\end{table}

\begin{figure}[h]
	\makebox[\textwidth][c]{
	\begin{subfigure}[t]{.6\textwidth}
		\centering
		\includegraphics[width=1.05\textwidth]{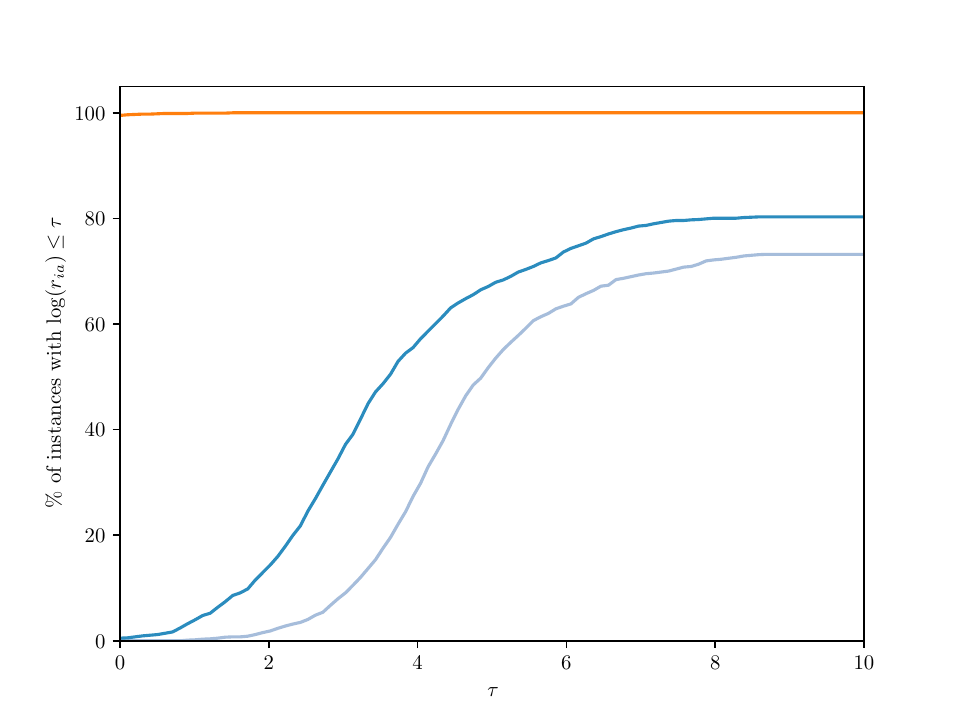}
		\caption{Performance profiles showing percentage of instances with a performance ratio of within $\tau$. Plotted on a log scale.}\label{fig:j10_performance_profile}
	\end{subfigure}%
	\qquad
	\begin{subfigure}[t]{.6\textwidth}
		\centering
		\includegraphics[width=1.05\textwidth]{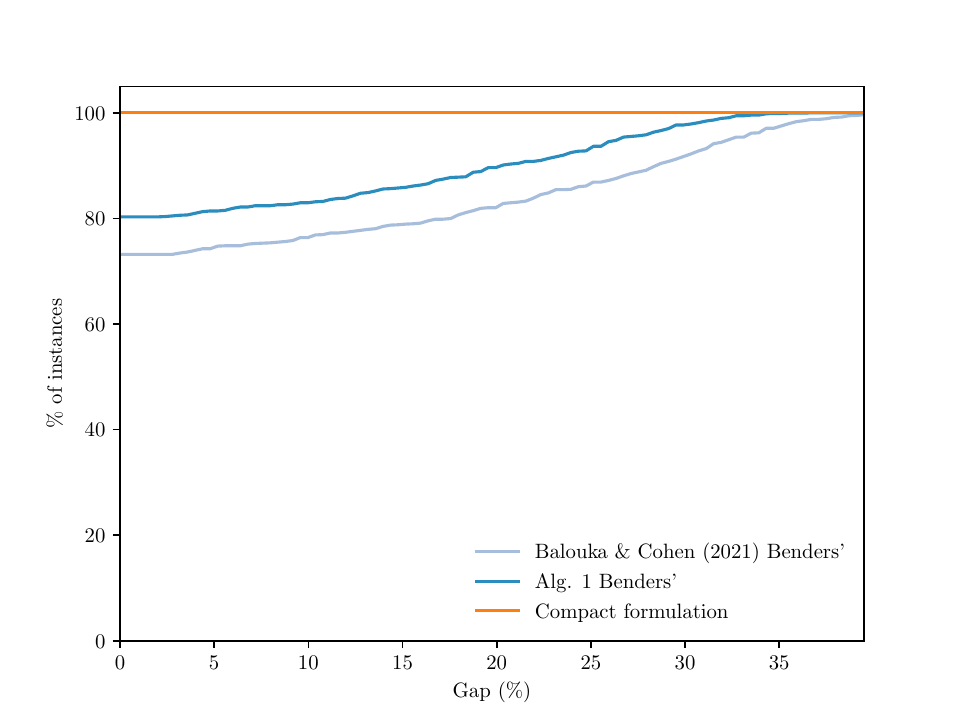}
		\caption{Percentage of instances solved to within given optimality gap within the two-hour time limit.}\label{fig:j10_gaps}
	\end{subfigure}
	}
	\caption{Comparison of the Benders' approaches and compact formulation over instances in the $j10$ set.}
	\label{fig:j10_plots}
\end{figure}

\begin{figure}[h]
	\makebox[\textwidth][c]{
	\begin{subfigure}[t]{.6\textwidth}
		\centering
		\includegraphics[width=1.05\textwidth]{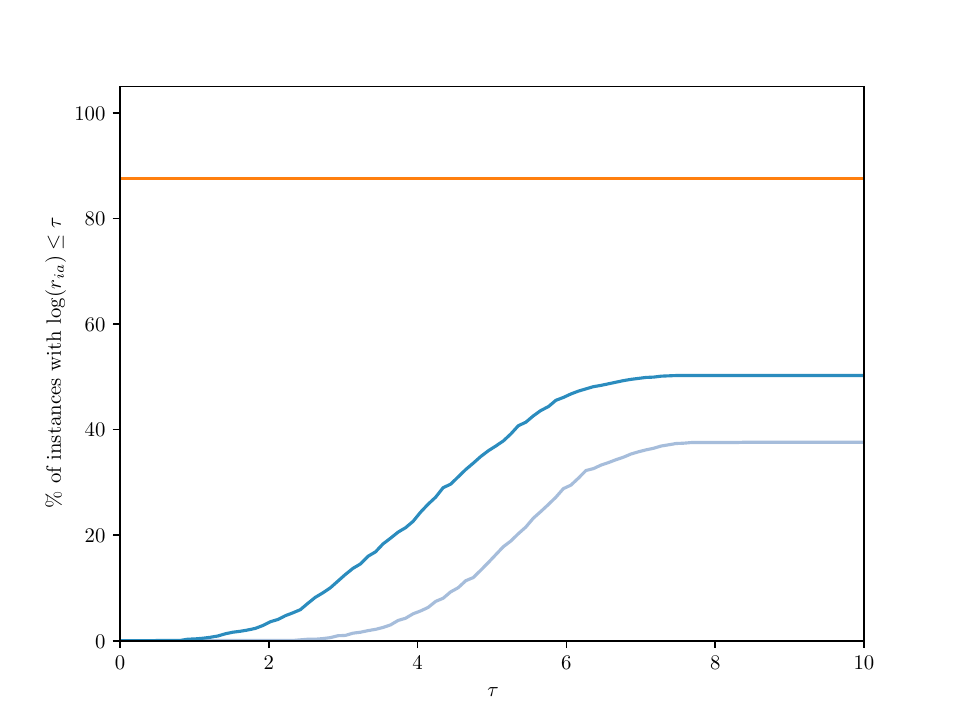}
		\caption{Performance profile showing percentage of instances with a performance ratio of within $\tau$. Plotted on a log scale.}\label{fig:j20_performance_profile}
	\end{subfigure}%
	\qquad
	\begin{subfigure}[t]{.6\textwidth}
		\centering
		\includegraphics[width=1.05\textwidth]{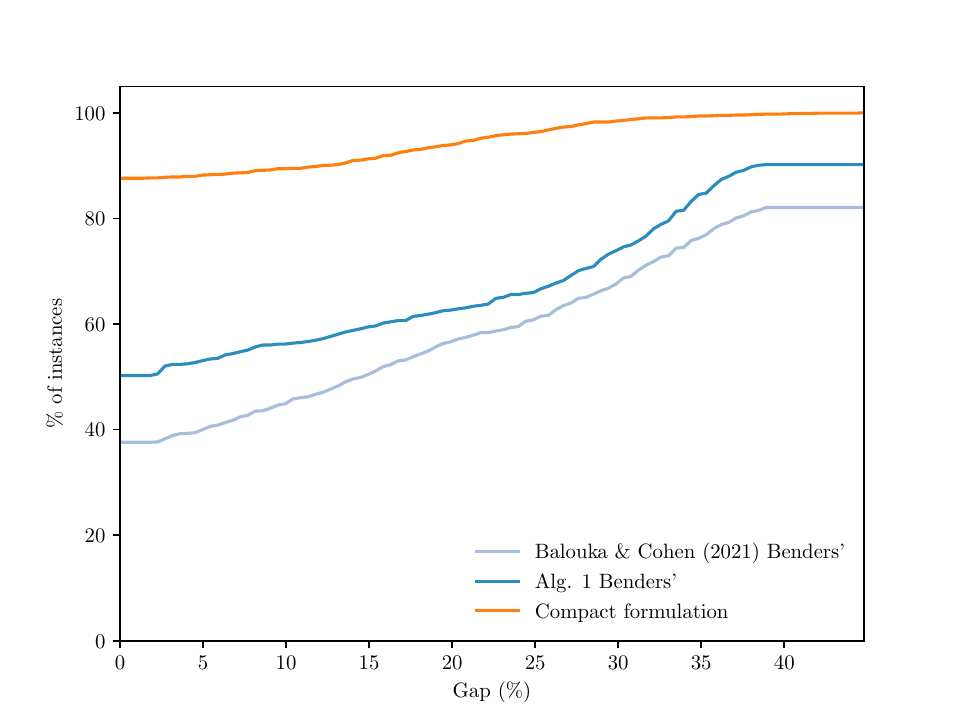}
		\caption{Percentage of instances solved to within given optimality gap within the two-hour time limit.}\label{fig:j20_gaps}
	\end{subfigure}
	}
	\caption{Comparison of the Benders' approaches and compact formulation over instances in the $j20$ set.}
	\label{fig:j20_plots}
\end{figure}

Figures \ref{fig:j10_plots} and \ref{fig:j20_plots} show performance profiles and optimality gap plots for the compact formulation and both implementations of the Benders' decomposition solution approach (Algorithm \ref{alg:benders} and \cite{balouka2021robust}), across the $j10$ and $j20$ sets respectively. Note that only the results for the instances with non-zero values for $\Gamma$ were used to generate these plots, resulting in a total of 1608 $j10$ instances and 1664 $j20$ instances. 

Performance profiles (\cite{dolan2002benchmarking}) provide a visual comparison of competing solution approaches using their \textit{performance ratios}. The performance ratio of algorithm $a\in \mathcal{A}$ for instance $i\in \mathcal{I}$ is defined to be 
$$r_{ia}=\frac{t_{ia}}{\min_{a\in \mathcal{A}} t_{ia}},$$
where $t_{ia}$ is the time required to solve instance $i$ using algorithm $a$. If method $a$ is unable to solve an instance $i$ within the 2-hour time limit, then $r_{ia}=R$, where $R\geq \max_{i,\,m}r_{im}$. The performance profiles in Figures \ref{fig:j10_performance_profile} and \ref{fig:j20_performance_profile} show the percentage of instances that have a performance ratio within a factor of $\tau$ of the best algorithm. These are plotted using a log scale for clarity. The y-intercepts of each method give the percentage of instances for which that method found the optimal solution in the fastest time, whereas the right-most value gives the overall percentage of instances that were solved to optimality within the time limit by that method.

These performance profiles show that whenever an instance was able to be solved to optimality by either solution approach, the compact formulation was always the faster method. If we specifically compare the solution times of the compact formulation with the stronger of the two Benders' implementations, across the instances in $j10$ for which both methods found the optimal solution, the average solution time of the compact formulation (1.4s) was almost 200 times faster than for the Benders' approach (283.5s). When both methods found the optimal solution for instances in the $j20$ set, the average solution time of the compact formulation (13.7s) was almost 100 times faster than for the Benders' approach (1304.5s). 

Figures \ref{fig:j10_gaps} and \ref{fig:j20_gaps} show the percentage of instances solved by each method to within a given optimality gap within the 2-hour time limit. These plots serve as a continuation of the performance profiles beyond just the instances that were solved to optimality. Summarising the data from these plots across both the $j10$ and $j20$ instance sets, the \cite{balouka2021robust} Benders' implementation solves 55.1\% of instances to optimality and finds feasible solutions for a further 35.6\% of instances. Our Benders' implementation solves 65.0\% of instances and finds feasible solution solutions for a further 30.0\% of instances. The compact formulation on the other hand finds a feasible solution to every instance, solving 93.1\% these to optimality.

The results presented here show the clear improvements to the Benders' algorithm afforded by amending the master problem to use the formulation (\ref{eqn:master_first})-(\ref{eqn:master_last}). More significantly however, these results demonstrate the complete dominance of the compact formulation over both Benders' methods.

\section{Conclusions}

The work presented in this paper extends the compact mixed-integer programming formulation first introduced by \cite{bold2021compact}, for application to the two-stage adjustable robust MRCPSP with uncertain activity durations. The computational performance of this formulation has been examined over a total of 3270 uncertain MRCPSP instances of varying size and difficulty, and compared against an improved version of the current state-of-the-art for solving this problem, based on a Benders' decomposition approach. The improved Benders' approach is the result of replacing the original master problem with a new formulation based on a simplified version of the compact formulation we present for the full problem. 

Results presented in Section \ref{section:experiments} show that the compact formulation completely dominates the enhanced Benders' approach, solving over 43\% more instances to optimality, and doing so with dramatically reduced computation times. In addition to these strong computational improvements, the proposed compact formulation has the significant added benefit of being simpler to implement than the iterative Benders' approach.

Despite these strong results, the instances used in these experiments contain only up to 20 activities. To enable the solving of larger scale instances, the development of heuristic solution approaches should be a primary focus of future research on this problem. 

\section*{Acknowledgements}

The authors are grateful for the support of the EPSRC-funded (EP/L015692/1) STOR-i Centre for Doctoral Training.

\bibliography{paper}

\end{document}